\documentclass[letterpaper, 10 pt, conference]{IEEEconf}
\IEEEoverridecommandlockouts

\usepackage[mathscr]{euscript}

 \usepackage{adjustbox}
\usepackage{amsmath,amssymb,amsfonts}
\usepackage[linesnumbered,ruled,vlined]{algorithm2e}
\SetKwInput{KwInput}{Input}                
\SetKwInput{KwOutput}{Output}  
\SetKwInput{kwProcedure}{Procedure}           
\usepackage{graphicx}
\usepackage{microtype}

\usepackage{hyperref}

\usepackage{multicol}

\usepackage[mathscr]{euscript}
\usepackage{bigints}
 
\usepackage{algpseudocode}
\usepackage{tcolorbox}

\usepackage{amssymb}
\makeatletter
\renewcommand*\env@matrix[1][\arraystretch]{%
  \edef\arraystretch{#1}%
  \hskip -\arraycolsep
  \let\@ifnextchar\new@ifnextchar
  \array{*\c@MaxMatrixCols c}}
\makeatother

\usepackage{bbm}

 \let\mathscr\relax
\usepackage[scr]{rsfso}
\usepackage{dirtytalk}
\usepackage{ntheorem}
\usepackage{tabularx}
\usepackage{hyperref}
\usepackage{outlines}
\usepackage{mathtools}
\usepackage{color}
\usepackage{float}
\usepackage{diagbox}
\usepackage[utf8]{inputenc}
\usepackage{lipsum}
\usepackage{graphicx}
\usepackage{subcaption}
\usepackage{cuted}
\usepackage{verbatim}
\newtheorem*{examp*}{Example}
\usepackage{marvosym}
\usepackage{array}
\usepackage{amsfonts}
\usepackage{enumerate}
\usepackage{graphicx}

\providecommand{\keywords}[1]{\textbf{\textit{Index terms---}} #1}

\newtheorem*{proof*}{Proof}
\newtheorem*{prob*}{Problem}

\newtheorem*{solution*}{Solution}
\newtheorem*{sol*}{Solution}
\newtheorem*{GOCP*}{GOCP}

\newtheorem{theorem}{Theorem}

\newtheorem{cop*}{Convex Optimization Problem}
\newtheorem*{ihlqr*}{IHLQR}
\newtheorem{definition}{Definition}
\newtheorem{proposition}{Proposition}
\usepackage{textcomp}

\begin{document}

\title{Choosing Augmentation Parameters in OSQP-\\
A New Approach based on Conjugate Directions\\
 
\thanks{\textit{Acknowledgement}: Part of this work was done in while the author was working at Télécom SudParis, Institut Polytechnique de Paris. The author is thankful to Prof. \href{https://nguyenhoainam.wp.imt.fr/}{Hoai Nam Nguyen} for the fruitful discussions thereat.}
}

\author{ Avinash Kumar\\ 
Email: avishimpu@gmail.com
}

\maketitle

\begin{abstract}
OSQP is a general purpose solver, based upon the alternating direction method of multipliers, for convex quadratic programs. Within this solver's algorithm, the idea of the augmented Lagrangian with a penalty parameter- a parameter which captures the relative weight-age on the objective function and the constraints of the problem in-hand- is utilized to develop an algorithm with so-called \textit{augmentation parameters}. The selection of these parameters is a crucial task and the optimal way to do the selection is not yet known. This work proposes a new method to select these parameters by utilising the information of the conjugate directions of the coefficient matrix of a linear system of equations present in the algorithm. This selection makes it possible to \textit{cache} these conjugate directions, instead of computing them at each iteration, resulting in a faster computation of the solution of the linear system thus reducing the overall computation time. This reduction is demonstrated by a numerical example by comparing the total time taken for the algorithms to converge sufficiently close to the optimal solution.
\end{abstract}

{\keywords  OSQP, Quadratic programs, ADMM, augmentation parameters, conjugate directions} 

\section{Introduction}

The analysis of convex optimization problems is of great importance in applied mathematics. Convex optimization problems arise in numerous fields including - machine learning, control engineering (model predictive control), Lasso and Hubber fitting and so on \cite{stellato2020osqp}. With the advent of big data and consequently higher dimensional data availability, the size of the convex optimization problems that are required to be solved continues to grow in size. Therefore, to handle real-time optimization-based tasks, it becomes necessary to develop algorithms that can work efficiently and effectively with the high-dimensional data. Furthermore, study of the solution algorithms for convex optimization problems turns out fruitful for the analysis of non-convex problems.

In this work, we consider a class of convex optimization problems called quadratic programs.
There exist a plethora of algorithms to solve the quadratic programs and can be broadly classified into three categories: $(1)$ active-set methods \cite{wolfe1959simplex, nocedal1999numerical}, $(2)$ interior-point methods \cite{nesterov1994interior}, and $(3)$ first-order methods \cite{frank1956algorithm}. The major drawback of active set methods is that their worst-case complexity grows exponentially with the number of constraints \cite{klee1972good} while the interior point methods are not scalable for very large scale problems because each iteration requires heavy computational tasks to be performed \cite{nguyen2022improved}.

In recent years ADMM, which is an operator splitting based first-order method,  has received a lot of attention by the researchers. This is because this method seems to be very well suited for the large scale problems because of the decomposability (inherited from the dual ascent algorithm) and  superior convergence properties (inherited from the method of multipliers) \cite{boyd2011distributed}. Moreover, the steps involved in the ADMM algorithm are computationally cheap and simple to implement \cite{stellato2020osqp}. However, the optimal selection of the parameters, inherited from the augmented Lagrangian, involved in the ADMM algorithm is still an open problem \cite{boyd2011distributed}. 

Operator Splitting Quadratic
Program (OSQP) is  general-purpose solver for quadratic programs developed upon the alternating direction method of multipliers (ADMM) as applicable to the convex optimization problems. Operator Splitting Quadratic Program (OSQP), proposed in \cite{stellato2020osqp}, is one of the most popular general-purpose solvers for the convex quadratic programs. The algorithm is based on a novel splitting technique and subsequent usage of the auxiliary variables leading to a quasi-definite linear system which is always solvable. 

One of the most computationally expensive step in the OSQP algorithm involves the solution of a linear system-a common trait of the algorithms derived from ADMM. It is this linear system which we propose to solve in a more efficient manner. To do so we use the information about the \textit{conjugate directions} of the coefficient matrix. The information about the conjugate directions is of grave importance in the study of convex optimization. To the best of author's knowledge, this line of research, where the information of the conjugate directions is exploited to improve the characteristics of the algorithm, remains unexplored in the literature and the author believes that it holds promise. The information about the conjugate directions allows us to compute offline a set of parameters and cache them thus avoiding their computation at each iteration of the algorithm leading to reduction in the computation time of the overall algorithm which is demonstrated by a numerical example. Note that in this study, we do not take into consideration the memory budget requirements and assume that sufficient storage capacity is available for the caching the computed augmentation parameters and conjugate directions.

The rest of the paper is organised as follows. Section II presents a review of the OSQP solver as reported in \cite{stellato2020osqp} followed by the background of conjugate direction and conjugate gradient methods. In Section III, we present the method to choose augmentation parameters by utilising the information of the conjugate directions of a specific coefficient matrix. Section IV presents the comparison of the proposed algorithm with the conjugate direction method. Section $V$ discusses the future aspects of the work.

 \textbf{Notation}:Throughout this manuscript, scalars and scalar-valued functions  are denoted by small ordinary alphabets, vectors are denoted by bold ordinary alphabets, and matrices are denoted by capital alphabets. $X \succ \mathbf{0}$ means that the matrix $X$ is symmetric positive definite. $diag(\bullet)$ represents a diagonal matrix. A `{$T$}' in the superscript represents the transpose. `{$^k$}' in the superscript of a quantity denotes the value of the quantity at the $k^{\text{th}}$ iteration.
\section{Background}
  
\subsection{OSQP}
OSQP is used for solving the quadratic programs as given by \eqref{COP}. 
\begin{equation}
\begin{gathered}
    \text{minimize }   0.5 \mathbf{x}^T P \mathbf{x} + \mathbf{q}^T \mathbf{x}   \\
\text{subject to }   \mathbf{l} \leq A \mathbf{x} \leq \mathbf{u}   
\label{COP}
\end{gathered}
 \end{equation}
where $\mathbf{x} \in \mathbb{R}^{n}$ is the decision variable, $P \in  \mathbb{R}^{n \times n}$ is a positive semi-definite matrix, $\mathbf{q} \in \mathbb{R}^{n }$, $A \in \mathbb{R}^{m \times n}$, and $\mathbf{l} \in  \mathbf{R}^{m}$ and  $\mathbf{u} \in \mathbb{R}^{m}$ define the lower and upper limits on the variable $A\mathbf{x}$. This formulation allows handling equality constraints by setting $\mathbf{l}_i= \mathbf{u}_i$. 
Let us define the set $\mathbb{C}$ as $\mathbb{C}= \{ \mathbf{z} \in \mathbb{R}^{m}| l_i \leq  z_i \leq u_i, i=1,2....., m \}$. With this definition and introduction of a new decision variable $\mathbf{z}$, we can rewrite the quadratic program \eqref{COP} in the form below which is more suitable for the application of a standard ADMM algorithm.
\begin{equation}
\begin{gathered}
        \text{minimize }  0.5 \mathbf{x}^T P \mathbf{x} + \mathbf{q}^T \mathbf{x} \\
    \text{subject to } A \mathbf{x}=\mathbf{z}; \mathbf{z} \in \mathbb{C}.
    \label{COP2}
\end{gathered}
\end{equation}    
 
The OSQP algorithm as applicable to the quadratic program \eqref{COP2} is detailed in Algorithm \ref{alg1}. Therein $\mathbf{y} \in \mathbb{R}^m$ is the Lagrange multiplier associated with the constraint $A \mathbf{x}= \mathbf{z}$, $\tilde{\mathbf{x}}$ and $\tilde{\mathbf{z}}$ are the auxiliary variables, $\sigma >0$ and  $\rho >0$ are the augmentation parameters, $\alpha \in (0,2)$ is the relaxation parameter and $\Pi$ denoted the Euclidean projection onto $\mathbb{C}$. More details of this algorithm can be found in \cite{stellato2020osqp}.

\begin{algorithm}

 \textbf{given}  initial values $\mathbf{x}^0$, $\mathbf{z}^0$, $\mathbf{y}^0$ and parameters $\rho>0$, $\sigma >0$, $\gamma \in (0,2)$; $k=0$
 
 \textbf{repeat}

$ (\tilde{\mathbf{x}}^{k+1}, \mathbf{\nu}^{k+1}) \leftarrow $
\text{(solve the linear system)} 
$$ \begin{bmatrix}P+ \sigma I & A^T \\
 A & - \rho^{-1} I 
 \end{bmatrix} \begin{bmatrix} \tilde{\mathbf{x}}^{k+1} \\ \mathbf{\nu}^{k+1} \end{bmatrix}
= \begin{bmatrix}\sigma \mathbf{x}^{k} - \mathbf{q} \\ \mathbf{z}^k - \rho^{-1} \mathbf{y}^k \end{bmatrix}
$$

$\mathbf{\tilde{z}}^{k+1} \leftarrow  \mathbf{z}^k+ \rho^{-1} \left( \mathbf{\nu}^{k+1} - \mathbf{y}^{k} \right)$

$ \mathbf{x}^{k+1} \leftarrow  \gamma  \tilde{\mathbf{x}}^{k+1} + \left( 1- \gamma  \right) \mathbf{x}^{k}$

$ \mathbf{z}^{k+1} \leftarrow \Pi \left( \gamma  \tilde{\mathbf{z}}^{k+1} + (1-\gamma ) \mathbf{z}^{k} + \rho^{-1} \mathbf{y}^k \right)$

$ \mathbf{y}^{k+1} \leftarrow \mathbf{y}^{k} + \rho \left( \gamma  \tilde{\mathbf{z}}^{k+1} +(1-\gamma ) \mathbf{z}^{k} - \mathbf{z}^{k+1}\right)$

\textbf{until} termination criterion is satisfied

\caption{OSQP Algorithm for QP \eqref{COP2}}
\label{alg1}
\end{algorithm}

It is well known that the most computationally expensive task in Algorithm \ref{alg1} is the 
step $3$ where we need to solve the linear system
\begin{equation}
    \begin{bmatrix}
        P+ \sigma I & A^T \\
        A & - \rho^{-1} I
    \end{bmatrix}
    \begin{bmatrix}
       \mathbf{ \tilde{x}}^{k+1} \\ 
        \mathbf{\nu}^{k+1}
    \end{bmatrix} = \begin{bmatrix}
    \sigma \mathbf{x}^k- \mathbf{q} \\
    \mathbf{z}^k - \rho^{-1} \mathbf{y}^k
    \end{bmatrix}.
    \label{eq_dir}
\end{equation}

How to choose the parameters $\sigma>0 $, $\gamma \in ( 0,2) $ and $\rho >0$ in the best possible way is still an open problem. The most crucial parameter is the augmentation parameter $\rho$ and it known that having different values of $\rho$ for different constraints can greatly improve the performance of the algorithm and thus it is reasonable to choose this augmentation parameter as a diagonal positive definite matrix instead of a scalar, let us denote this matrix as $\varrho \in \mathbb{R}^{m \times m}$, $\varrho =diag(\rho_1, \rho_2, \hdots \hdots, \rho_m)$ with $\rho_i>0 \ \forall i \ \in \{1,2,3, \cdots, m\}$ which we call the augmentation matrix in the rest of the manuscript. Furthermore, it is beneficial to have an \textit{adaptive} $\varrho$ which updates at each iteration. We shall revisit this point when we propose a way to compute the elements of $\varrho$. Moreover, for large scale quadratic programs, we can use the following equations  instead of equation \eqref{eq_dir}, obtained by eliminating $\mathbf{\nu}^{k+1}$ and utilization of the Schur complement lemma. This approach is called the indirect method \cite{stellato2020osqp}.
 
\begin{equation}
\left(P+ \sigma I +  A^T \varrho A \right) \tilde{\mathbf{x}}^{k+1}=\sigma \mathbf{x}^k- \mathbf{q} + A^T \left( \rho \mathbf{z}^k - \mathbf{y}^k\right). \label{eq_inv}
\end{equation}
\begin{equation*}
\tilde{\mathbf{z}}^{k+1}= A \mathbf{\tilde{x}}^{k+1}.
\end{equation*}
In this work, we propose a method to efficiently solve \eqref{eq_inv} using the information of the conjugate directions of the positive-definite matrix $P+ \sigma I +  A^T \varrho A$.

 \subsection{Conjugate Direction Methods}
 In this section, we provide a background of conjugate directions methods as applicable to the area of optimization \cite{nocedal1999numerical}.

\begin{definition}[Conjugate Directions] {\cite{nocedal1999numerical}}
    The set of nonzero vectors $\{ \mathbf{p}_1, \mathbf{p}_2, \hdots \hdots , \mathbf{p}_n \}$, $\mathbf{p}_i \in \mathbb{R}^{n}$, is said to be conjugate with respect to a symmetric positive definite matrix $\bar{A} \in \mathbb{R}^{n \times n}$ if $\mathbf{p}_i^T\bar{A} \mathbf{p}_j=0$ $\forall$ $i \neq j$. The vectors $\{ \mathbf{p}_1, \mathbf{p}_2, \hdots \hdots , \mathbf{p}_n \}$ are called the conjugate directions of the matrix $\bar{A}$.
    \label{def_1}
\end{definition}

 The method of conjugate directions as applicable for obtaining a solution to the system 
 \begin{equation}
      \bar{A}\bar{\mathbf{x}}=\bar{\mathbf{b}} \text{ where }  \bar{A} \succ \mathbf{0} \text{ and } \bar{\mathbf{x}} \in \mathbb{R}^n.
      \label{eq_lin}
 \end{equation} 
 is detailed in Algorithm \ref{alg2}. The solution to \eqref{eq_lin} is $\bar{\mathbf{x}}^*= \bar{A}^{-1} \bar{\mathbf{b}}$.

\begin{algorithm}
    \textbf{given} initial value $\mathbf{\bar{x}}_0 \in \mathbb{R}^{n}$ and the set of conjugate directions $\{ \mathbf{p}_1, \mathbf{p}_2, \hdots \hdots , \mathbf{p}_n \}$; $k=0$
    
    \textbf{repeat}
    $$\bar{\alpha}_k= - \frac{\left( \bar{A} \bar{\mathbf{x}}_k - \bar{\mathbf{b}} \right)^T \mathbf{p}_k}{\mathbf{p}_k^T \bar{A} \mathbf{p}_k}$$  

   $$  \bar{\mathbf{x}}_{k+1}= \bar{\mathbf{x}}_k + \bar{\alpha}_k \mathbf{p}_k$$

   \textbf{until} the termination criteria is satisfied
   \caption{Conjugate Directions Method for \eqref{eq_lin}}
   \label{alg2}
\end{algorithm}

\begin{theorem} \cite{nocedal1999numerical}
    For any $\mathbf{\bar{x}}_0 \in \mathbb{R}^n$ the sequence ${\bar{\mathbf{x}}_k}$ generated by the conjugate direction algorithm converges to the solution $\bar{\mathbf{x}}^*$ of the linear system in at most $n$ steps.
    \label{thm_cd}
\end{theorem}

The importance of conjugate directions is clear from Theorem \ref{thm_cd}- the solution $\bar{\mathbf{x}}^*$ can be arrived at in at-most $n$ iterations if the steps are taken along the conjugate directions $\{ \mathbf{p}_1, \mathbf{p}_2, \hdots \hdots , \mathbf{p}_n \}$ of $\bar{A}$. But the computation of the complete set of conjugate directions requires an excessive amount of computation. Considering this point and the fact that there is no best known method to choose augmentation matrix $\varrho \in \mathbb{R}^{m \times m}$  in \eqref{eq_inv} for utilisation in Algorithm \ref{alg1}, in this work, we propose \textit{a way to choose the augmentation matrix $\varrho$ in such a way that conjugate directions of coefficient matrix $P+\sigma I + A^T \varrho A$ can be computed offline and cached}. We detail the method of making such a selection in the next section.

A more popular algorithm closely related to the conjugate directions algorithm is the conjugate gradient method in which the conjugate directions are computed while progressing through the iterations using \textit{only} the previous direction and thus the storage of all the directions is not required hence reducing the computational burden and storage requirements \cite{nocedal1999numerical}. In this work, we use the version of the conjugate gradient method as detailed in Algorithm \ref{alg3} and as applicable for solving the linear system \eqref{eq_lin} \cite{nocedal1999numerical}.

\begin{algorithm}
    \textbf{given} initial value $\mathbf{\bar{x}}_0 \in \mathbb{R}^{n}$ , \textbf{set} $\bar{\mathbf{r}}_0 \leftarrow \bar{A} \bar{\mathbf{x}}_0 - \bar{\mathbf{b}}$, $\bar{\mathbf{p}}_0\leftarrow -\bar{\mathbf{r}}_0$; $k=0$
    \textbf{repeat}    
    $ \bar{\alpha}_k \leftarrow \dfrac{\bar{\mathbf{r}}_k^T \bar{\mathbf{r}}_k}{\bar{\mathbf{p}}_k^T \bar{A} \bar{\mathbf{p}_k}}$
  
   $ \bar{\mathbf{x}}_{k+1} \leftarrow \bar{\mathbf{x}}_k+ \bar{\alpha}_k \bar{\mathbf{p}}_k$
   
   $ \bar{\mathbf{r}}_{k+1} \leftarrow \bar{\mathbf{r}}_k+ \bar{\alpha}_k \bar{A} \bar{\mathbf{p}}_k$
   
   $\bar{\beta}_{k+1} \leftarrow \dfrac{\bar{\mathbf{r}}_{k+1}^T \bar{\mathbf{r}}_{k+1}}{ \bar{\mathbf{r}}_k^T \bar{\mathbf{r}}_k}$

   $\bar{\mathbf{p}}_{k+1} \leftarrow  -\bar{\mathbf{r}}_{k+1} + \bar{\beta}_{k+1} \bar{\mathbf{p}}_k$

   $k \leftarrow k+1$
 
   \textbf{until} the termination criteria is satisfied    (here $\bar{\mathbf{r}}_k= \bar{A} \bar{\mathbf{x}}_k- \bar{\mathbf{b}}$ is the residual at $k^{th}$ iteration and $\bar{\mathbf{p}}_k$s are the conjugate directions.)
   
   \caption{Conjugate Gradient Method for \eqref{eq_lin}}

   \label{alg3}
\end{algorithm}

\section{Main Results}

Proposition \ref{prop1} helps choose the parameter $\varrho$ such that the conjugate directions of $P+ \sigma I +  A^T \varrho A$ can be stored. 

Let $\{ \mathbf{e}_1, \mathbf{e}_2, \hdots \hdots \mathbf{e}_n \}$, $\mathbf{e}_i \in \mathbb{R}^{n}$, be the set of standard basis vectors of $\mathbb{R}^{n}$. Let the $\{ \mathbf{d}_1, \mathbf{d}_2, \mathbf{d}_3, \hdots \hdots \mathbf{d}_n \}$, $\mathbf{d}_i \in \mathbb{R}^n$, be a set of conjugate directions of the matrix $ \left(P+ \sigma I +  A^T \varrho A \right)$. Let $\mathbf{a}_i$ denote the $i^{th}$ column of matrix $A^T$. Decomposing these quantities, we get $ \mathbf{a}_i= \sum_{k=1}^{n} a_{ik} \mathbf{e}_{k}$, and $d_p=\sum_{k=1}^{n} \alpha_k^p \mathbf{e}_k$ where $a_{ik}$s and $\alpha_{k}^p$s are scalars.
\begin{proposition}
\label{prop1}
The conjugate directions $\mathbf{d}_p \left(=\sum_{k=1}^{n} \alpha_k^p e_k \right)$ and $\mathbf{d}_q \left(=\sum_{k=1}^{n} \alpha_k^q e_k \right)$ ($p \neq q$ $(p,q) \in \{1,2,3 \hdots \hdots, n \}$) of the matrix $(P+ \sigma I +  A^T \varrho A)$ satisfy \eqref{eq_1}.
 \end{proposition}
\begin{equation}
\begin{gathered}
\left(\sum_{k=1}^{n}\alpha^p_k e_k^T \right) \left( P+ \sigma I \right) \left( \sum_{k=1}^{n} \alpha_k^q e_k \right)+  \\
\sum_{l=1}^{m} \rho_l\left( \sum_{k=1}^{n} \alpha_k^p a_{lk}^2 \alpha_k^{q}+  \sum_{m,j=1;i \neq j}^{n} a_{lm}a_{lj} \alpha_m^{p}\alpha_j^{q}
\right)=0.
\end{gathered}
\label{eq_1}
\end{equation}
 
\begin{proof*}
\normalfont
The result follows by Definition \ref{def_1} and subsequent simplifications.
$\square$
\end{proof*}

\subsection{Choosing the set of Conjugate Directions}
Solving equations \eqref{eq_1} gives a set of conjugate directions are well as the augmentation parameters $\rho_i$s, $i \in \{1,2,3, \cdots ,m \}$. These quantities can then be stored and utilized in the online phase of the algorithm. However, solving \eqref{eq_1} a fixed (constant) value of the augmentation parameters is obtained which may not be desirable as varying the parameter $\varrho$ with iterations helps in improving the performance of the algorithm. In this work, we utilize a slightly different version of the equations as stated in \eqref{eq_p_si}.
 
\begin{equation}
\begin{gathered}
\left(\sum_{k=1}^{n}\alpha^p_k e_k^T \right) \left(P+ \sigma I\right) \left( \sum_{k=1}^{n} \alpha_k^q e_k \right)=0; \\
\sum_{l=1}^{m} \rho_l\left( \sum_{k=1}^{n} \alpha_k^p a_{lk}^2 \alpha_k^{q}+  \sum_{m,j=1;i \neq j}^{n} a_{lm}a_{lj} \alpha_m^{p}\alpha_j^{q}
 \right)=0.
 \label{eq_p_si}
\end{gathered}
\end{equation}
 
Clearly, if $ \rho_i$s satisfy \eqref{eq_p_si} then $a \rho_is$ also satisfy \eqref{eq_p_si} for a scalar $a \in \mathbb{R}$- we call this feature the \textit{scalability-invariance} in terms of $\rho_i$s of the solutions of \eqref{eq_p_si}. 
Also, it is easy to see that \eqref{eq_p_si} is special case of \eqref{eq_1} and hence solutions to \eqref{eq_p_si} verify the desired properties as possessed by solutions of \eqref{eq_1}. The offline phase of the proposed algorithm is as summarized below. It is worth noting that the selection of $\varrho$ and the conjugate directions is now not independent of each other.
\begin{tcolorbox}
\textbf{Offline Computation}

Solve the nonlinear algebraic equations  \eqref{eq_p_si} to obtain $\varrho$, $\alpha^p$s and $\alpha^q$s, $p \neq q$, $p,q \in \{1,2,3, \hdots , n \}$ and cache the conjugate directions $\mathbf{d}_i= \sum_{k=1}^{n} \alpha_k^i \mathbf{e}_k $, $i \in \{1,2,3, \hdots, n \}$ and the augmentation matrix say $\varrho =\varrho_{off} \succ \mathbf{0}$.
\end{tcolorbox}

\subsection{Adaptive Update of $\varrho$ }
Updating $ \varrho$ at each iteration offers faster convergence rates of the algorithm. Considering this fact, and observing that the equations \eqref{eq_p_si} hold even if $\varrho$ is scaled at later iterations, we utilise an update rule for $\varrho$  proposed in \cite{boyd2011distributed} which is based upon the primal and dual residuals at the current iteration. The initialisation and update of the $\varrho$ matrix as proposed in \cite{stellato2020osqp} is as detailed below.

 \textbf{Initialisation:} $\varrho= diag(\rho_1, \rho_2 \hdots, \rho_m )$, 
\begin{align}
\rho_i= \begin{cases}
    \bar{\rho}   \text{ if }  l_i \neq u_i \\
    10^3 \bar{\rho}   \text{ if } l_i=u_i,
\end{cases}
\label{rho_init_boyd}
\end{align}
Here $\bar{\rho}>0$.

\textbf{Update-scheme:}
\begin{equation}
\begin{adjustbox}{max width=350pt}
$
    {\varrho}^{k+1} \leftarrow {\varrho}^{k}  \sqrt{\frac{\vert| \mathbf{r}^{k}_{prim} \vert|_{\infty}/ \max \{ ||{A} \mathbf{x}^k \vert|_{\infty}, \vert|| \mathbf{z}^k \vert| |_{\infty}  \}} {|| \mathbf{r}^k_{dual}||_{\infty} / \max\{||{P} \mathbf{x}^k ||_{\infty}, 
||{A}^T \mathbf{y}^k ||_{\infty}, || \mathbf{q}||_{\infty}\}
}}
$.
\end{adjustbox}
\label{eq_upd_aug}
\end{equation}
where $\mathbf{r}^k_{prim}$ and $\mathbf{r}^k_{dual}$ are the primal and dual 
at iteration $k$ and are as given below.
$$
   \mathbf{r}^k_{prim}= A \mathbf{x}^k- \mathbf{z}^k,\\
    \mathbf{r}^k_{dual}= P \mathbf{x}^k + \mathbf{q}+A^T \mathbf{y}^k.
$$
where $\mathbf{y} \in \mathbb{R}^{m}$ is the Lagrange multiplier associated with the constraint $A \mathbf{x}= \mathbf{z}$. 

Contrariwise, in this work we use $\varrho_{off}$ as computed in the offline phase to initialize $\varrho$. Then, thanks to the \textit{scalability-invariance} feature of the solutions of \eqref{eq_p_si}, the scheme \eqref{eq_upd_aug} can be adopted to update the augmentation matrix $\varrho$ at each iteration and \eqref{eq_p_si} still holds. It is also known that if the $\varrho$ is updated at each iteration but the updates stop after a fixed number of iterations, then the convergence can be proved. Thus, there is one more free parameter in the algorithm-say $N_c>2$- the number of iterations after which the updates on augmentation parameter stop. A rigorous analysis of this parameter is one of the topics of the future work. After the offline computation, we have obtained a set of the conjugate directions of the matrix $P+ \sigma I +  A^T \varrho A$ and the corresponding augmentation parameter matrix $\varrho$. Once these quantities are computed, the following algorithm is run online to obtain the solution of the optimization problem at each iteration. We now use this information in the standard OSQP algorithm to invert the matrix $P+ \sigma I +  A^T \varrho A$  using the conjugate directions method where the conjugate directions have been cached. This drastically improves the computation time taken for inverting the matrix.

\begin{tcolorbox}
\textbf{Online Computation}
\begin{enumerate}
    \item 
\textbf{given:} initial values $\mathbf{x}_0$, $\mathbf{z}_0$, $\mathbf{y}_0$, $\sigma >0$, $\alpha \in (0,2)$, $N_c>2$
the set of conjugate directions $\mathbf{d}_i$s $i \in \{1,2,3,....,n \}$ and the augmentation matrix $\varrho= \varrho_{off} \succ \mathbf{0} \in \mathbb{R}^{m \times m}$ as computed in the offline phase
    \item[] \textbf{ repeat}
     \item Solve the system 
    $$(P+ \sigma I+  A^T \varrho A) \tilde{x}^{k+1}= \sigma x^{k} -q +A^T \left( \rho z^k -y^k\right)$$ using Algorithm \ref{alg2} by utilising the conjugate directions $\mathbf{d}_k$s as computed in the offline phase.
        $\tilde{z}^{k+1}\leftarrow A \tilde{x}^{k+1}$.
   \item    if $iteration\_count<N_c$, update $\varrho$ as per \eqref{eq_upd_aug}.
   \item[]   \textbf{ until} the termination criterion is satisfied perform the steps $5$, $6$ and $7$ in Algorithm \ref{alg1}.
\end{enumerate}
\end{tcolorbox}

\section{Numerical Example}
In this section, we present the comparison of the proposed algorithm with the conventional OSQP algorithm by solving \eqref{eq_inv} utilising 
\begin{enumerate}
    \item the conjugate gradient method i.e. Algorithm \ref{alg3}, and
    \item the proposed approach.
\end{enumerate}

Consider the quadratic program \eqref{COP} with the matrices
with $
    P=
\begin{bmatrix}
    3  &   1    & 3 &    2 \\
     1  &   1   &  2 &    1 \\
     3   &  2  &   8  &   4\\
     2    & 1 &    4   &  3
\end{bmatrix}  $, $A= I$, $\mathbf{q}=\begin{bmatrix}  1 & 1 & 1 & 1  \end{bmatrix}^T$, $\mathbf{l}= \begin{bmatrix} -2 & -1 & -3 & -4 \end{bmatrix}^T$ and $\mathbf{u} =\begin{bmatrix} 10 & 1 & 3 & 0\end{bmatrix}^T$.

\begin{figure}
    \includegraphics[scale=0.17]{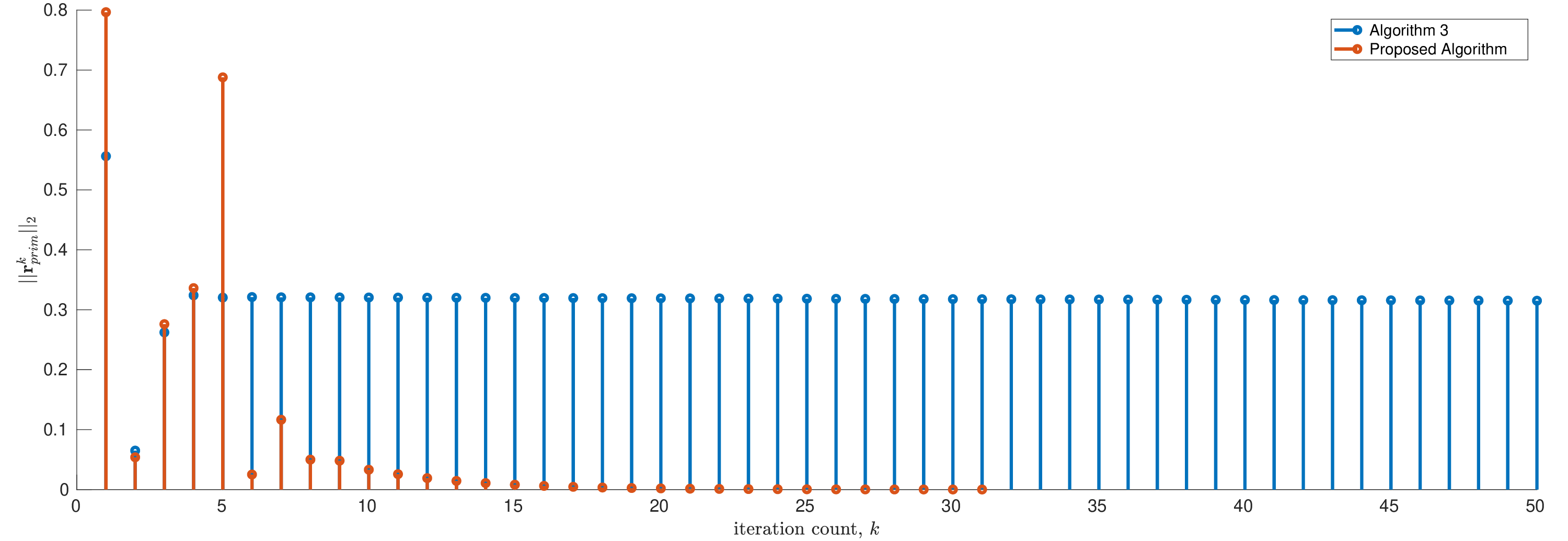}
    \caption{Variation of the norm of the primal residual}
    \label{fig1}
\end{figure}

\begin{figure}[H]
    \includegraphics[scale=0.17]{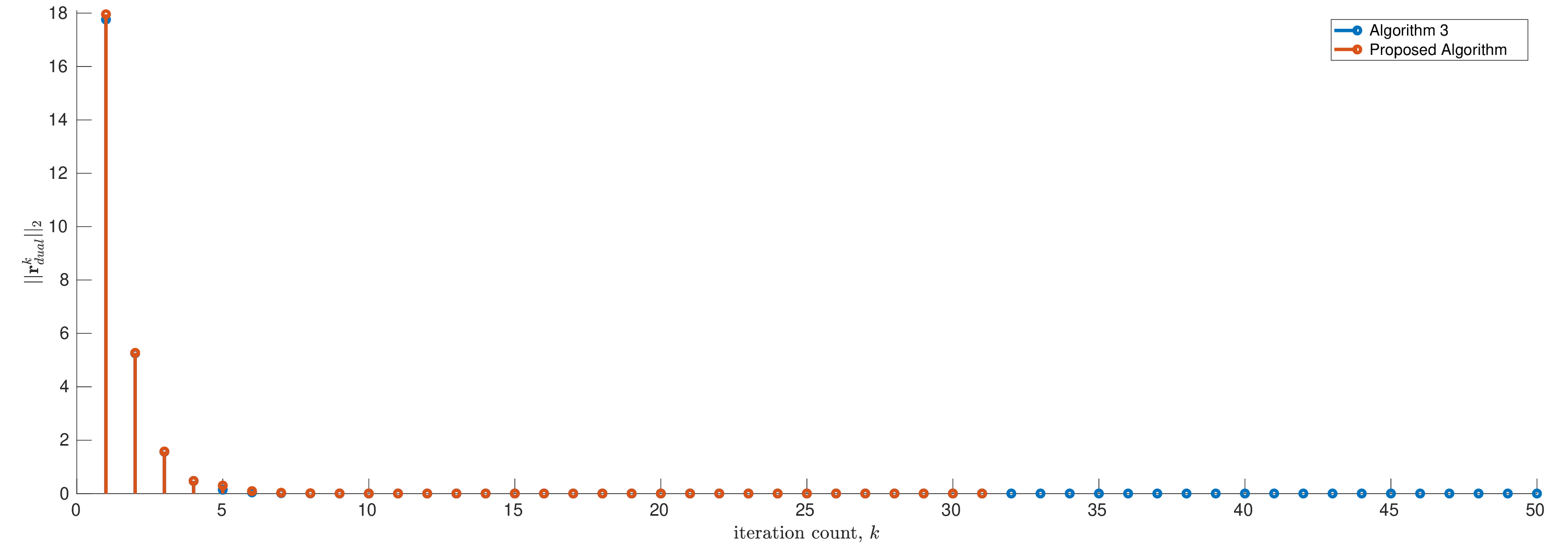}
    \caption{Variation of the norm of the dual residual}
    \label{fig2}
\end{figure}

\begin{enumerate}
    \item \textbf{When the augmentation matrix $\varrho$ is initialised as $\varrho_{off}$. }
    
\textit{Proposed Method:}
The offline computation of the proposed for this problem gives the following results.
\begin{itemize}
    \item Conjugate Directions
    
$
        d_1=\begin{bmatrix}
              0.2605&
    2.5954&
    0.4708&
   -2.3267
        \end{bmatrix}^T$,$d_2= \begin{bmatrix}
                0.3084 &  0.1520 &    0.3328&    0.2067
        \end{bmatrix}^T$, $d_3= \begin{bmatrix}
                2.4128&  -0.1005&   -1.1522&   -0.1601 
        \end{bmatrix}^T$ and
        $d_4= \begin{bmatrix}
               -0.3493&    1.1900 &   -0.5754 &    0.7348
        \end{bmatrix}^T$.
    \item Augmentation matrix, 
    $$\varrho_{off}= 
    \begin{bmatrix}
        0.1000 & 0 & 0& 0\\
        0 & 0.1087 & 0 & 0 \\
        0 & 0& 0.1757&  0\\
         0& 0& 0& 0.1631\end{bmatrix}.$$
\end{itemize}

We set $N_c=5$, $\sigma=0.0001$ and $\gamma=1.3$. The termination criteria used is $|| \mathbf{r}_{dual}||_2< 10^{-4}$ and $|| \mathbf{r}_{prim}||_2< 10^{-4}$.

We compare the results with the conjugate gradient based approach where $\varrho$ is initialised with $\varrho_{off}$.
To compare the results, we average the three parameters- $(1)$ total time to reach the solution ($T_{tot}$), $(2)$ time taken to solve the linear system involved ($T_{inv}$) and $(3)$ the number of iterations $N$ taken to obtain the solution obtained using the two algorithms over $10000$ runs with the initial conditions $\mathbf{x}_0$ derived from a pseudo-normal distribution by utilising the MATLAB function $randn$. It is easily seen that the $T_{inv}$ is lesser for the proposed approach which leads to lesser overall computation time $T_{tot}$.

\item \textbf{When the augmentation matrix is initialised as per \eqref{rho_init_boyd} with $\bar{\rho}=0.2$.}

For the case when the augmentation parameter is initialized with \eqref{rho_init_boyd} for the conjugate gradient based approach, the variation of primal and dual residual norms for the two algorithms are as shown in Figure \ref{fig1}  and Figure \ref{fig2} respectively, for an initial condition $ \mathbf{x}_0=\begin{bmatrix} 1 &2& 3& 4 \end{bmatrix}^T$. It is easily seen that the residuals vanish faster for the proposed approach where $\varrho$ is initialised with $\varrho_{off}$.

\end{enumerate}

   \begin{table}
   \vspace{0.3 cm}
\begin{center}
    \begin{tabular}{|c|c|c|c|}
    \hline
    \textbf{Technique Used} &$T_{tot} (ms)$&$T_{inv}(ms)$ & $N$\\
    \hline
        Algorithm \ref{alg3}
        ($\varrho$ initialised with $\varrho_{off}$) &     $0.3236$  &  $0.0036$ &  $ 33$    \\ 
    \hline
    Proposed Algorithm  &  $0.2913$ &    $0.0026$   &   $33$     \\ 
    \hline
\multicolumn{4}{|c|}{ms=milliseconds} \\
\hline
    \end{tabular}
    \caption{Comparison of different approaches}
    \label{tab:my_label}
    \end{center}
\end{table}

(The simulations were done on MATLAB R2024a on a system with processor 11th Gen Intel® Core™ i5-1135G7 @ 2.40GHz × 8 and 16GB RAM.)
\section{Conclusion}
In this work, we present a new method to choose the augmentation matrix $\varrho$ in the OSQP algorithm, this selection is based upon the idea of the conjugate directions. The selection is made in such a way that the conjugate directions of a specific matrix involved in the linear system of equations can be cached thus facilitating solving the linear system involved in the iterations, thus reducing the overall computational complexity of the algorithm. The numerical example demonstrates the efficacy of the proposed approach in terms of the reduction in the computation time. The future work includes more rigorous treatment of the proposed algorithm along with the development of a suitable solver.
\bibliographystyle{unsrt}

\end{document}